\documentstyle[12pt,amstex,amssymb]{article}
\title{Optimal Prediction and the Klein-Gordon Equation}
\author{{\sc Ole H. Hald}\bigskip\\Department of
Mathematics\\University of California\\
Berkeley, California 94720--3840\bigskip\\{\tt
hald@@math.berkeley.edu}\\phone 510--642--4809\\fax:
510--642--8204\bigskip\bigskip\bigskip\\
{\footnotesize Communicated by Alexandre J. Chorin,}\\
{\footnotesize Department of Mathematics, University of
California, Berkeley, CA 94720--3840}}
\date{}
\newtheorem{lemma}{Lemma}
\newtheorem{theorem}{Theorem}
\begin{document}
\setlength{\baselineskip}{22pt}
\maketitle

\vfill{}
\noindent Classification: Applied Mathematics

\noindent Number pages: 16

\noindent Number of words in abstract: 58

\noindent Total number characters in paper:
$\sim$44,000

\newpage

\begin{abstract}
The method of optimal prediction is applied to calculate the future means of
solutions to the Klein-Gordon equation. It is shown that in an appropriate
probability space, the difference between the average of all solutions that
satisfy certain constraints at time $t=0$, and the average computed by
an approximate method, is small with high probability.
\end{abstract}

\newpage

\setcounter{page}{1}
\section{Introduction}
The method of optimal prediction was introduced by Chorin, Kast,
Kupferman \cite{CK1,CK2,CK3} to study complicated flows,
hopefully including turbulence at a future
time.
Instead of solving a particular initial value problem we ask for
the average of all solutions that satisfy certain constraints at
time $t=0$.
The constraints may be local averages of the initial data, or a small
number of Fourier coefficients.
Neither will determine the initial data uniquely.
The idea then is to use statistical information to compensate for the
incompleteness of the initial data.
In its most elementary version the method of optimal prediction is more
expensive than solving the original initial value problem.
The savings are achieved by finding an evolution equation for the
constraints and from this determining the average of the solutions
for $t > 0$.
For non-linear problems this can only be done approximately.
However, for linear problems we can estimate the difference between
the exact averages and the averages computed by the approximate method.
We get the sharpest bound if the constraints are close to an invariant
subspace for the adjoint of the differential equation.
We apply the theory to the Klein-Gordon equation and prove that the difference
between the exact mean at time $t$ and the outcome of an approximate
calculation is small with high probability.
We also show that the exact averages converge with probability $1$ as
we increase the dimension of the trial space.
This remains true even if the measure is carried by weak solutions that
are difficult to obtain individually.
We confine ourselves to a single case, but the arguments can be extended
to the linear Schr{\"o}dinger equation and to linear Korteveg de Vries
equations.
\section{Two Methods}
In this section, we will present an exact and an approximate method for
finding the average of the solutions to a differential equation.
Let $L$ be a real $m \times m$ matrix and let $G$ be a real $m \times n $
matrix of rank $ n < m $.
We will look at the solutions $u(t)$ of
\begin{equation}
\label{1}
{\dot u (t) = L u(t)}
\end{equation}
and assume that the initial conditions satisfy the constraint
\begin{equation}
\label{2}
{G^T\! u(0) = v_0 .}
\end{equation}
If $S(t) = e^{tL} $ is our fundamental matrix, then $u(t) = S(t) u(0) $.
To find the average of all $u$ that satisfy (\ref{2})
we need a measure. Let $A$ be a positive definite
matrix of order $n$ and define
\begin{displaymath}
{ P( u \in B ) = \int_B Z^{-1} \ e^{-\frac{1}{2} u^T\! A u } \ du , }
\end{displaymath}
where $Z$ is chosen so that $ P ( {\Bbb R}^m ) = 1 $.
If $L^T\! A + AL = 0 ,$ then $P$ is an invariant measure, i.e.
$P(B)=P(S(t)B)$ for all $t$.
The matrix $A$ may be chosen in many ways, but there is a natural choice
if (\ref{1}) is a Hamiltonian system.
By restricting $P$ to the set $G^T\! u = v_0$ and normalizing again, we get
a measure $P'$ that satisfies
\begin{equation}
\label{3}
{\langle u \rangle = \int_{G^T\! u = v_0} u \ dP' = A^{-1} G M^{-1} v_0 ,}
\end{equation}
where $M = G^T\! A^{-1} G $, see \cite {CK1,CK2,CK3}.
Since $u(t) = S(t) u(0) $ we can determine the average of all solutions
that satisfy $G^T\! u(0) = v_0$ and get
\begin{equation}
\label{4}
{ {\langle u(t) \rangle}_{exact} = S(t) \langle u(0) \rangle =
S(t) A^{-1} G M^{-1} v_0 }
\end{equation}

The approximate method is harder to motivate.
We would not expect that $G^T\! u(t) = v_0$ for all $t>0$;
but there may exist a function $v(t)$ such that
$G^T\! u(t) = v(t) $
for all $u(t)$ that satisfy $G^T\! u(0) = v_0$.
The arguments for $t=0$ are then applicable.
After replacing $v_0$ in (\ref{3}) by $v(t)$ we see
that
$\langle u(t) \rangle = A^{-1} G M^{-1} v(t) $.
In addition,
$v(t) = G^T \langle u(t) \rangle $,
and it follows from (\ref{1}) that
$\dot v(t) = G^T\! L \langle u(t) \rangle $.
We can now formulate the approximate method. Let $K=LA^{-1}$. Then
\begin{equation}
\label{5}
{{\langle u(t) \rangle}_{approx} = A^{-1} G M^{-1} v(t) }
\end{equation}
\begin{equation}
\label{6}
{\dot v(t) = G^T\! K G M^{-1} v(t), \qquad v(0) = v_0 \ .}
\end{equation}

If $n \ll m$, it should be cheaper to find the approximate solution than
the exact solution. The question is: "How good is the approximation?".
To answer this question, we set
\begin{displaymath}
{e(t) = {\langle u(t) \rangle}_{approx} -  {\langle u(t) \rangle}_{exact}}
\end{displaymath}
\begin{displaymath}
{ E = L^T\! G + G M^{-1} G^T\! K G \, . }
\end{displaymath}
Suppose $L^T\!A+AL=0$. Then $A^{-1}L^T\!+LA^{-1}=0$, and it follows from
(\ref{4}), (\ref{5}), (\ref{6}) that
\begin{eqnarray*}
\dot e (t) \!\!&=& \!\!A^{-1}GM^{-1} \dot v(t) - \dot S(t) A^{-1} G M^{-1} v_0
\\
&=& \!\! A^{-1}GM^{-1}G^T\!KGM^{-1}v(t) - L S(t) A^{-1}GM^{-1} v_0
\\
&=& \!\! A^{-1} G M^{-1} G^T\! K G M^{-1} v(t) -
 L [ A^{-1} G M^{-1} v(t) - e (t) ]
\\
&=& \!\!Le(t) + A^{-1} [ L^T\! G+G M^{-1} G^T\! KG ] M^{-1} v(t) \ .
\end{eqnarray*}
Using the explicit solution of inhomogeneous linear equations, (see \cite{CL}
page 78), we obtain
\begin{equation}
\label{7}
e(t) = \int_0^t S(t-s) A^{-1} E M^{-1} v(s) ds .
\end{equation}

\begin{lemma}
\label{lemma1}
If $L^T\! A + AL = 0 $, then
\begin{displaymath}
{ | A^{1/2} e(t) | \leq t \ | A^{-{1/2}} E M^{-1/2} |
\ | M^{-1/2} v_0 | .}
\end{displaymath}
\end{lemma}
{\bf Proof:} To bound $e(t)$, we need two facts:
\begin{equation}
\label{8}
{(A^{1/2}S(t)A^{-1/2})}^T (A^{1/2}S(t)A^{-1/2}) = I
\end{equation}
\begin{equation}
\label{9}
v^T\! (t) M^{-1} v(t) = v_0^T \! M^{-1} v_0.
\end{equation}
Equation (\ref{8}) says that $A^{1/2} S(t) A^{-1/2} $
is orthonormal, while (\ref{9}) corresponds to
conservation of energy for (\ref{6}). Both are
consequences of the assumption $L^T\!A+AL=0$. To prove
(\ref{8}), we differentiate with respect to $t$, use
$\dot S = LS $, and obtain
\begin{displaymath}
\frac{d}{dt} [ A^{-1/2} S^T\! (t) A S(t) A^{-1/2} ] =
A^{-1/2} S^T\! (t) [ L^T\!A+AL ] S(t) A^{-1/2} = 0 .
\end{displaymath}
The matrix $A^{-1/2} S^T\! (t) A S(t) A^{-1/2} $ is therefore independent
of time and is equal to the identity when $t=0$.
To prove (\ref{9}), we differentiate with respect to
$t$, use (\ref{6}) and $K^T\!+K=0$, and get
\begin{displaymath}
\frac{d}{dt} [ v^T\! (t) M^{-1} v(t) ] =
v^T\! (t) M^{-1} G^T (K^T\! + K ) G M^{-1} v(t) = 0.
\end{displaymath}
This shows that $v^T\!(t) M^{-1} v(t) $ is independent of time.
We can now complete the proof of Lemma \ref{lemma1}.
Multiplying both sides of (\ref{7}) by $A^{1/2}$ and
using
(\ref{8}), (\ref{9}) yield
\begin{eqnarray*}
|A^{1/2} e(t) | & \leq & \int_0^t
|A^{1/2} S(t-s) A^{-1/2} | \ | A^{-1/2} E M^{-1/2} | \
| M^{-1/2} v(s)|\ ds \\
& \leq & t\ |A^{-1/2} E M^{-1/2} | \ | M^{-1/2} v_0 |.
\end{eqnarray*}
This completes the proof.

It follows from Lemma \ref{lemma1} that $e(t) \equiv 0$
if $E=0$.
This will occur if $G$ is a left invariant subspace for $L$.
To prove this, let $L^T \!G = GB $.
Then $G^T\! A^{-1} L^T\! G = G^T\! A^{-1} G B $, and we see that
$B = - M^{-1} G^T\! K G $ and $E = L^T\!G+GM^{-1}G^TKG = L^T\!G-GB = 0$.

\section{ Hamiltonian Systems}
It is not true that for every $L$ there is a positive definite matrix $A$
such that $L^T\! A + A L = 0 $.
You need the eigenvalues of $L$ to be purely imaginary
and that $L$ be diagonalizable.
However, $A$ exists for linear Hamiltonian systems.
Lets look at $\ddot q (t)= - A_0^2 q(t) $,
where $A_0$ is positive definite.
This equation describes small oscillations around equilibrium.
Setting $\dot q(t)=p(t)$, we arrive at
\begin{displaymath}
{ \frac{d}{dt} \left[ \!\!\begin{array}{c} q(t) \\ p(t) \end{array}
\!\!\right] =
\left[ \begin{array}{cc} 0 & I \\ {-A^2_0} & 0 \end{array}
\right]
\left[ \!\!\begin{array}{c} q(t) \\ p(t) \end{array}\!\!\right].}
\end{displaymath}
The Hamiltonian for this system is
$h = \frac{1}{2} [ p^T\!p+q^T\! A^2_0 q ] $,
i.e. $\dot q_i = \partial_{p_i} h $ and
$\dot p_i = - \partial_{q_i} h $.
It is natural to constrain $p$, $q$ separately
\begin{displaymath}
\left[ \begin{array}{cc} G^T_q & 0 \\ 0 & G^T_p \end{array}\right]
\left[ \!\!\begin{array}{c} q(0) \\ p(0) \end{array}\!\!\right] =
\left[ \!\!\begin{array}{c} v_q(0) \\ v_p(0) \end{array}\!\!\right] .
\end{displaymath}
More complicated relations between $p(0)$, $q(0)$ are possible and may
be preferable in special cases.
Letting
$ u(t) = { q(t) \brack p(t) } $,
we have $\dot u = Lu$, $G^T\!u(0) = v_0$, and $h= \frac{1}{2} u^T\! A u$
as in (\ref{1}), (\ref{2}) where
\begin{displaymath}
L = \left[ \begin{array}{cc} \ & I \\ -A_0^2 & \ \end{array}\right] , \
G = \left[ \begin{array}{cc} G_q & \ \\ \ & G_p \end{array}\right] , \
A = \left[ \begin{array}{cc} A^2_0 & \ \\ \ & I \end{array}\right] .
\end{displaymath}

Set $|u|_A = |A^{1/2} u | = (2h)^{1/2} $.
Since $M = G^T\! A^{-1} G $ and $K = LA^{-1}$, we obtain
\begin{displaymath}
M = \left[ \begin{array}{cc} G^T_q\! A_0^{-2} G_q & \ \\ \ &
G^T_p \!G_p \end{array}\right], \;\;\;
G^T\! KG = \left[ \begin{array}{cc} \ & G^T_q\! G_p \\ -G^T_p\!
G_q  & \
\end{array}\right] .
\end{displaymath}
Note that $M$ is positive definite and that $G^T\!KG$ is skew symmetric.
To simplify the analysis, we assume that $G_p = G_q = G $ and hope
that the double use of $G$ will not cause confusion.
The differential equation for the approximate method can then be written
as
\begin{equation}
\label{10}
{ \frac{d}{dt} \left[\!\! \begin{array}{c} v_q(t) \\ v_p(t) \end{array}
\!\!\right] =
\left[ \begin{array}{cc} \ & I \\ {-(G^T\! G)(G^T\!A_0^{-2}G)^{-1}} & \
\end{array}
\right]
\left[ \!\!\begin{array}{c} v_q(t) \\ v_p(t) \end{array}\!\!\right] , }
\end{equation}
$\mathrm{cf.}$ (\ref{6}).
If $G$ consists of eigenvectors of $A_0^2$, then each
eigenfrequency of (\ref{10}) agree with an
eigenfrequency of the original problem and $e(t)=0$.
To estimate the error in the approximate method, we must bound
$|A^{-1/2} E M^{-1/2} |$ in {Lemma \ref{lemma1}}.
Since $ E = L^T\! G + G M^{-1} G^T\! K G $, it follows that
$ A^{-1/2} E M^{-1/2} =
{\: 0 \: F \brack 0 \ 0 } $,
where
\begin{equation}
\label{11}
F = -A_0 G (G^T\!G)^{-1/2} + A_0^{-1} G (G^T\! A_0^{-2} G )^{-1}
(G^T\!G)^{1/2}.
\end{equation}
Thus, $|A^{-1/2} E M^{-1/2} |= |F|$, and it is enough to bound the 2-norm
of
\begin{equation}
\label{12}
F^T\!F = (G^T\!G)^{-1/2} (G^T\! A_0^2 G) (G^T\!G)^{-1/2} -
(G^T\!G)^{1/2}(G^T\! A_0^{-2} G )^{-1} (G^T\! G)^{1/2} .
\end{equation}
To continue the analysis, we turn to a specific problem.

\section{Klein-Gordon}
In the paper by Chorin, Kast, Kupferman \cite{CK1,CK2,CK3}
the method of optimal
prediction was applied to linear and non-linear Schr{\"o}dinger
equations.
Here we will study the Klein-Gordon equation
\begin{equation}
\label{13}
u_{tt} = u_{xx} - u
\end{equation}
on the interval $0 \leq x \leq 2 \pi $ with periodic boundary conditions.
The equation describes dispersive waves on a string subject to a restoring
force.
A similar equation occurs in relativistic quantum field theory \cite{LHR}.
The Hamiltonian for (\ref{13}) is
\begin{equation}
\label{14}
h(t) = \frac{1}{2} \int_0^{2 \pi} (u_t)^2 + (u_x)^2+(u)^2 \ dx .
\end{equation}
The corresponding Hamiltonian system is
\begin{equation}
\label{15}
\partial_t \left[ \!\!\begin{array}{c} u(x,t) \\ \pi (x,t) \end{array}
\!\!\right] =
\left[ \begin{array}{cc} 0 & I \\ \partial_x^2-I & 0 \end{array} \right]
\left[ \!\!\begin{array}{c} u(x,t) \\ \pi (x,t) \end{array} \!\!\right] ,
\end{equation}
where $\pi (x,t) = u_t (x,t) $.
Note that $A_0^2 = -\partial_x^2+I$.
For a derivation see \cite{Gol}.
 We constrain the initial data by
prescribing local averages around the points
$ x_{\alpha} = 2 \pi \alpha / (2n+1) $ for $\alpha = 0,1, \ldots ,2n. $
Specifically,
\begin{eqnarray}
\label{16}
\int_0^{2 \pi} g(x-x_{\alpha} ) u(x,0) dx &=& v_{q, \alpha } (0), \\
\int_0^{2 \pi} g(x-x_{\alpha} ) \pi(x,0) dx &=& v_{p, \alpha } (0). \nonumber
\end{eqnarray}

Let us imagine that $v_p(0)$, $v_q(0)$ are given, and set
$v_0 = { v_q(0) \brack v_p(0) } $.
Following Chorin, Kast, Kupferman \cite{CK1,CK2,CK3}, we let
\begin{equation}
\label{17}
g(x) = \frac {1}{\sqrt{2 \pi }} \sum_{k=-\infty}^{\infty}
e^{-k^2\! \sigma^2 \!/4} \frac{e^{i k x }}{\sqrt{2 \pi}}.
\end{equation}
The function $g$ is positive, and $2 \pi$ periodic, has norm $1$  and decrease
away from the origin.
As $\sigma \rightarrow 0 $, $g$ tends to a delta function.
Since the measure $P$ is finite dimensional, we assume that there is an
integer $m \geq 0$ such that all $u(x,t)$, $\pi (x,t)$ can be written as
\begin{displaymath}
\sum_{k=-m}^{m} c_k  \frac{e^{i k x }}{\sqrt{2 \pi}},
\end{displaymath}
where $\bar{c}_k = c_{-k}$ and $m = n + r(2n+1)$.
The complex notation is equivalent to
\begin{displaymath}
\frac{a_0}{\sqrt{2 \pi}} +
\sum_{k=1}^{m} \left( a_k \frac{\cos kx }{\sqrt{\pi}} +
b_k \frac{\sin kx}{\sqrt{\pi}} \right)
\end{displaymath}
when $c_0=a_0$ and $c_k=(a_k-i b_k)/\sqrt{2}$ for $k=1,2, \ldots ,m$.
In the expansion of $g(x)$, we replace $\exp (-k^2 \sigma^2\! /4 )$ in
(\ref{17}) by $0$ if $|k| > m $, thus obtaining $Proj_m
\ g(x)$. Our basic variables are not the trigonometric
functions, but their Fourier coefficients.
Let $(a_i ,b_i )$ be the Fourier coefficients for $u(x,t)$, and let
$(\alpha_i , \beta_i )$ be the Fourier coefficients for $\pi (x,t) $.
Set $q^T = (a_{m}, \ldots ,a_0,b_1 , \ldots , b_m ) $
and $p^T =(\alpha_m , \ldots , \alpha_0 , \beta_1, \ldots ,\beta_m )$.
We can then rewrite (\ref{15}) as
\begin{equation}
\label{18}
\frac{d}{dt} \left[ \!\!\begin{array}{c} q(t) \\ p(t) \end{array} \!\!\right] =
\left[ \begin{array}{cc} 0 & I \\ {- \Lambda^2} & 0 \end{array} \right]
\left[ \!\!\begin{array}{c} q(t) \\ p(t) \end{array} \!\!\right] ,
\end{equation}
where $ \Lambda = \mathrm{diag} \ (\omega_{m} , \ldots,
\omega_0, \ldots , \omega_m ) $
and $\omega _k^2 = k^2+1$.
Observe the shift in notation: the constants $m,n$ from Section 2 have
been replaced by $2(2m+1)$, $2(2n+1)$.
To find the analogue of (\ref{16}), we expand $u(x,t)
$ in a complex Fourier series, use (\ref{17}), and get
\begin{displaymath}
v_{q,\alpha}(0) = \frac{1}{\sqrt{2 \pi}} \sum_{\ell=-m}^{m}
e^{-\ell^2 \! \sigma^2\!/4} \ e^{ i \ell x_\alpha } \ c_\ell .
\end{displaymath}
Since the points $x_\alpha$ are equidistant we have an aliasing effect.
Let $\ell = k + {j(2n+1)}$ with $-n \leq k \leq n $ and $ -r \leq j \leq r $.
Then
\begin{eqnarray*}
v_{q,\alpha}(0) & = & \sum_{k=-n}^n \frac{e^{i k x_\alpha}}{\sqrt{2n+1}}
\cdot \sqrt{ \frac{2n+1}{2 \pi} } \
\sum_{j=-r}^r e^{-[k+j(2n+1)]^2\!\sigma^2\!/4} \ c_{k+j(2n+1)}
\\
& = & \sum_{k=-n}^n U_{\alpha k} \cdot w_k .
\end{eqnarray*}
Note that $\bar{w}_k = w_{-k}$.
The matrix $U$ is the building block for the discrete Fourier transform
and is unitary.
Set
\begin{displaymath}
\Gamma = \sqrt{ \frac{2n+1}{2\pi}} \ \mathrm{diag} \left(
e^{-m^2\!\sigma^2\!/4}, \ \ldots \ , 1, \ \ldots \ ,
e^{-m^2\!\sigma^2\!/4}
\right).
\end{displaymath}

If $ {c}^T = ( c_{-m}, \ldots, c_0, \ldots, c_{m} ) $,
we can write (\ref{16}) as
$v_q(0) = U [ I \ \cdots \ I ] \Gamma {c} $
with $2r+1$ blocks of $I$'s.
To express the constraints as a product of real matrices we let
$X$, $Y$ be of order $2n+1$
and $2 m +1 $, respectively, and of the form
\begin{displaymath}
\frac{1}{\sqrt 2}
\left[
\begin{array}{ccccccc}
1& \ & \ & \ & \ & \ & i \\
\ & \ 1 & \ & \ & \ & i & \ \\
\ & \ & 1 & \ & i & \ & \ \\
\ & \ & \ & { \sqrt{2}} & \ & \ & \ \\
\ & \ & 1 & \ & -i & \ & \ \\
\ & 1 & \ & \ & \ & -i & \ \\
1 & \ & \ & \ & \ & \ & -i
\end{array}
\right] .
\end{displaymath}
Note that $X$, $Y$ are unitary.
The matrix $Q=UX$ is orthonormal and the $\alpha$'th row of $Q$ is
\begin{displaymath}
\sqrt{\frac{2}{2n+1}}
\left[ \cos (nx_\alpha ), \ \ldots \ , \cos (x_ \alpha ) ,\
\frac{1}{\sqrt{2}} ,
\ \sin (x_\alpha ) , \ \ldots \ , \sin (nx_\alpha) \right].
\end{displaymath}
Since ${c} = Y q $ and $\Gamma Y = Y \Gamma $, we finally obtain
$v_q(0) = UX X^\ast [ I \ \cdots \ I ] Y \Gamma q = Q Z^T \Gamma q $.
Because $v_q(0)$, $Q$, $\Gamma$, $q$ are real, $Z$ must also be real.
Let $G^T = Q Z^T \Gamma = U [ I \ \cdots \ I ] Y \Gamma $.
The analogue of (\ref{16}) is then
\begin{equation}
\label{19}
\left[
\begin{array}{cc}
G^T & \ \\
\ & G^T
\end{array}
\right]
\left[ \!\!
\begin{array}{c}
q(0) \\ p(0)
\end{array}
\!\!\right] =
\left[\!\!
\begin{array}{c}
v_q(0) \\ v_p(0)
\end{array}
\!\!\right] .
\end{equation}

We can now solve (\ref{18}), (\ref{19}) by the exact
method (\ref{4}) and by the approximate method
(\ref{5}). To estimate the difference, we use Lemma
\ref{lemma1} and need the  following result.
\begin{lemma}
\label{lemma2}
If $n \geq 1$ and $(2n+1)\sigma^2 \geq 2 $, then
\begin{displaymath}
|A^{-1/2} E M^{-1/2} | \leq (1.6) \, (2n+1) \, e^{-(2n+1)\sigma^2\!/4}
\end{displaymath}
\end{lemma}
{\bf Proof:} To bound $|F|$, we will determine $G^T\!G$,
$G^T\!A_0^2 G$, $G^T\!A_0^{-2}G$ in (\ref{12})
explicitly. Observe that $A_0 = \Lambda $.
By using the complex representation of $G$, $Y \Gamma = \Gamma Y$,
and $Y Y^\ast = I $, we see that $G^T\!G = U D_1U^\ast$ where
\begin{displaymath}
D_1 = \frac{2n+1}{2\pi} \mathop{ \mathrm{diag} } \limits_{-n \leq k \leq n}
\left( \sum_{j=-r}^r e^{-[k+j(2n+1)]^2\!\sigma^2\!/2}\right).
\end{displaymath}
Interchanging $k$, $j$ with $-k$, $-j$ shows that $(D_1)_{-k} = (D_1)_k$
which implies that $X^\ast D_1 = D_1 X^\ast$.
Since $U=QX^\ast$, we conclude that $G^T\!G = QD_1Q^T $.
Similar arguments give $G^T\!A_0^2 G = Q D_2 Q^T$ and
$G^T\!A_0^{-2} G = Q D_3 Q^T$, where
\begin{displaymath}
D_2 = \frac{2n+1}{2\pi} \mathop{ \mathrm{diag} } \limits_{-n \leq k \leq n}
\left( \sum_{j=-r}^r e^{-[k+j(2n+1)]^2\!\sigma^2\!/2}
\{[k+j(2n+1)]^2+1\} \right)
\end{displaymath}
\begin{displaymath}
D_3 = \frac{2n+1}{2\pi} \mathop{ \mathrm{diag} } \limits_{-n \leq k \leq n}
\left( \sum_{j=-r}^r e^{-[k+j(2n+1)]^2\!\sigma^2\!/2}
\{[k+j(2n+1)]^2+1\}^{-1}\right).
\end{displaymath}
We can now determine (\ref{12}) explicitly.
Since $Q$ is orthonormal, it follows that
\begin{displaymath}
F^T\!F = Q [ D_1^{-1} D_2 - D_1 D_3^{-1} ] Q^T.
\end{displaymath}
If $r=0$ then $E=F=0$ and the approximate and exact method agree.
Let $r \geq 1 $ and suppose that the largest term in the diagonal matrix
$ D_1^{-1} D_2 - D_1 D_3^{-1} $ occurs in the $k$'th position.
Set $d_j = \exp (-[k+j(2n+1)]^2\!\sigma^2\!/2) $ and
$\lambda_j = [k+j(2n+1)]^2 +1 $.
Extracting the leading order term in each sum, we get
\begin{eqnarray*}
( D_1^{-1} D_2 - D_1 D_3^{-1})_{k} & = &
\frac{\sum d_j \lambda_j}{\sum d_j} -
\frac{\sum d_j}{\sum d_j \lambda_j^{-1}} \\
& = & \frac{ (d_0 \lambda_0 + a )(d_0 \lambda_0^{-1} + b ) - (d_0+c)^2 }
{(d_0+c) \ (d_0 \lambda_0^{-1}+b) } \\
& = & \frac {a(d_0 \lambda_0^{-1}+b) - d_0(c-\lambda_0 b ) - c ( d_0+c)}
{(d_0 \lambda_0^{-1}+b) \ (d_0+c) } .
\end{eqnarray*}
Since $c > \lambda_0 b $ and the $2$-norm is invariant under orthonormal
transformations, we see that
\begin{equation}
\label{20}
|F|^2 \leq \frac{a}{d_0} = \sum_{|j|=1}^r
e^{\{k^2-[k+j(2n+1)]^2\} \! \sigma^2 \!/2 }
\{[k+j(2n+1)]^2+1\} .
\end{equation}
To estimate the exponential, we observe that
\begin{eqnarray}
\label{21}
\lefteqn{k^2-[k+j(2n+1)]^2 \leq -j^2(2n+1)^2+2|k||j|(2n+1) \nonumber } \\
& = & -(j^2-|j|)(2n+1)^2-|j|(2n+1)-2|j|(n-|k|)(2n+1) .
\end{eqnarray}
Combining (\ref{20}), (\ref{21}) with $|k| \leq n $
results in
\begin{displaymath}
|F|^2 \leq \sum_{j=1}^r e^{-(2n+1)\sigma^2\!/2}
e^{-(j^2-j)(2n+1)^2\sigma^2\!/2}
\, 2 \,[j^2(2n+1)^2+k^2+1].
\end{displaymath}
Since $k^2+1 \leq (2n+1)^2/4$ when $|k| \leq n $ and $n \geq 1 $, we
conclude that
\begin{displaymath}
|F|^2 \leq (2n+1)^2 e^{-(2n+1)\sigma^2\!/2}
\sum_{j=1}^r e^{-(j^2-j)(2n+1)^2\!\sigma^2\!/2} \, 2 \, (j^2+1/4).
\end{displaymath}
The last sum is less than $2.53$ when $(2n+1)\sigma^2\!/2 \geq 1$ and
$n \geq 1 $. This completes the proof.

\section{Stochastic Convergence}
By combining Lemma \ref{lemma1} and Lemma \ref{lemma2} we can bound
the difference between the exact and the approximate method.
Since $|M^{-1/2}v_0|$ depends on $2n+1$, $\sigma$, and $v_0$, we have not
established convergence.
Suppose $v_p(0)$, $v_q(0)$ are generated by two particular
random functions $u$, $\pi $ with $u(x,0)$ looking like
Brownian motion, and $\pi (x,0) $ resembling white noise.
We can then show that the approximate method is close to the exact method
if $n$ is large.
The rate of convergence is high if there
is a substantial overlap of the kernels in the constraints.
To measure the error we use the norm $| \cdot |_A$ whose square equals
twice the total energy.
\begin{theorem}
\label{theorem1}
Let $n \geq 1 $, and assume that $(2n+1)\sigma^2 \geq 6(\nu+1) \log(2n+1) $
with $\nu>0$.
Let $p$, $q$ be picked at random with respect to $P$, and set $v_p(0)=G^T\!p$,
$v_q(0)=G^T \! q$.
Consider all solutions of (\ref{15}) that satisfy
(\ref{16}). Then
\begin{displaymath}
\left| {\left \langle \!\!\left(
\!\!\!\begin{array}{c} u(x,t) \\ \pi (x,t) \end{array}
\!\!\!\right) \!\!\right \rangle }_{exact} -  \;\
{\left \langle \!\!\left(
\!\!\!\begin{array}{c} u(x,t) \\ \pi (x,t) \end{array}
\!\!\!\right) \!\!\right \rangle }_{approx} \right|_{A} \leq
\frac{2.3 \ t}{(2n+1)^\nu}
\end{displaymath}
with probability greater than $1-(2n+1)^{-\nu} $.
\end{theorem}
{\bf Proof: } It follows from Lemmas \ref{lemma1} and
\ref{lemma2}
that
\begin{displaymath}
|A^{1/2} e(t) | \leq (1.6) \, t \, (2n+1) \, e^{-(2n+1)\sigma^2\!/4}
|M^{-1/2} v_0 |,
\end{displaymath}
where $ v_0 = {v_q(0) \brack v_p(0)} $.
To complete the proof, we use Chebyshev's inequality.
Let $E$ be the expected value corresponding to $P$. Since $2.3 > 1.6 \sqrt{2}$
and $\sigma$ is bounded below, we obtain
\begin{eqnarray}
P \left( |A^{1/2} e(t) | > \frac{2.3 \ t }{(2n+1)^\nu} \right) &\leq&
P \left( |M^{-1/2} v_0 | > \frac{\sqrt{2} \ e^{(2n+1)\sigma^2\!/4}}
{(2n+1)^{\nu+1} }\
\right) \nonumber\\
&\leq& \frac{E(|M^{-1/2}v_0|^2)}{2(2n+1)^{\nu+1}} .
\label{22}
\end{eqnarray}
Using the definition of $M$ from Section 3 in conjunction with $A_0 = \Lambda $
and (\ref{19}), we get
\begin{displaymath}
E(v_0^T\!M^{-1}v_0) = E[ (\Lambda q )^T\! \Lambda^{-1} G
(G^T\! \Lambda^{-2} G )^{-1} G^T\! \Lambda^{-1} (\Lambda q ) ]
+ E[ p^T\!G (G^T\!G)^{-1} G^T\!p ].
\end{displaymath}

Since $\Lambda $ is diagonal, the measure $P$ is given by
\begin{eqnarray*}
dP &=& Z^{-1} e^{-\frac{1}{2}[a_0^2+\alpha_0^2+\sum_{k=1}^m
\omega_k^2(a_k^2+b_k^2) + (\alpha_k^2+\beta_k^2) ]}
\ da_0 \cdots d\beta_m \\
Z &=& 2 \pi \prod_{k=1}^m \left( \frac{2\pi}{\omega_k^2} \cdot 2\pi\right).
\end{eqnarray*}
The components of $\Lambda q $ and $p$ are therefore independent Gaussian
random variables with mean 0 and variance 1, and it follows that
\begin{eqnarray*}
E(v_0^T\!M^{-1}v_0)&=& \mathrm{tr} \
[ \Lambda^{-1} G ( G^T\!\Lambda^{-2} G )^{-1}
G^T\! \Lambda^{-1} ] +
\mathrm{tr} \ [ G(G^T\!G)^{-1}G^T] \\
&=& \mathrm{tr} \ [ (G^T\!\Lambda^{-2}G)^{-1} (G^T\!\Lambda^{-2}G)]
+ \mathrm{tr} \ [ (G^T\!G)^{-1} (G^T\!G) ] \\
&=& 2(2n+1)
\end{eqnarray*}
Here $\mathrm{tr} = \mathrm{trace} $ and we have used
$\mathrm{tr} \ (AB) = \mathrm{tr}\ ( BA) $ if $A$ is an $n \times m $
matrix and
$B$ is $m \times n$.
Combining the last result with (\ref{22}) and taking
the complementary event finishes the proof.

We remark that the components of $v_0$ are strongly correlated.
Indeed, it follows from (\ref{19}) that
\begin{displaymath}
E(v_0v^T_0) =
\left[ \begin{array}{cc} G^T\!\Lambda^{-1} & \ \\
\ & G^T \end{array} \right]
E
\left[ \begin{array}{cc} (\Lambda q )(\Lambda q )^T  & \ \\
\ & p p^T \end{array} \right]
\left[ \begin{array}{cc} \Lambda^{-1} G & \ \\
\ & G \end{array} \right]
= M.
\end{displaymath}
Using the spectral decomposition of $G^T\!G$, we can calculate
the variances explicitly and get
\begin{eqnarray*}
\mathrm{var} [v_{p,\alpha}(0) ] &=&
\sum_{k=-n}^n (Q_{\alpha,k} )^2 ( D_1)_k
= \frac{1}{2\pi} \sum_{\ell =-m}^{m}
e^{-\ell^2\!\sigma^2\!/2} \\
&=& \int_0^{2\pi} [ Proj_m \ g(x) ]^2 \ dx.
\end{eqnarray*}
The variance of $v_{p,\alpha}(0)$ is therefore of order
$1/(\sqrt{2\pi} \sigma)$.
For $v_{q,\alpha}(0)$ we get an additional factor of
$\{\ell^2+1\}^{-1}$, and
$1/(2\pi) < \mathrm{var} \ (v_{q,\alpha} (0) ) < \coth ( \pi )/2 $.

Suppose the components of $v_0$ are chosen as independent,
normally distributed random variables with mean $0$ and variance $1$.
If $n\geq 4$ and $(2n+1)\sigma^2 \geq 6(\nu+1)
\log (2n+1)$, we can show that
any interval longer than $4$ contains points $t$ for which
\begin{displaymath}
E \left( |e(t) |^2_A \right) \geq (2n+1)^{(\nu+1)(2n+1)/4-1} .
\end{displaymath}
The initial constraint $v_0$ must therefore be consistent with the
mathematical model if we want convergence.

\section{Convergence in $L^2$}
In Section 5, we compared the outcome of two numerical methods.
Both are defined on finite dimensional spaces and involve a
finite number of Fourier coefficients.
What happens if we fix the number of constraints, but increase
the dimension of the space?
Each random choice of the Fourier coefficients
${ \{ a_i, b_i, \alpha_i, \beta_i \} }_{i=0}^{\infty}$
yields a sequence of constraint values.
Such a sequence may or may not converge.
We will show that the sequence of exact solutions generated by the
constraints converges with probability $1$.
Note we are not comparing results for different values of $n$.
They differ by large amounts.

Let $m = n+r(2n+1) $.
By solving (\ref{18}), (\ref{19}) explicitly and using
(\ref{4}), we  find that the Fourier coefficients for
the average of all solutions of
(\ref{13}) with the constraints (\ref{14}), satisfy
\begin{eqnarray*}
\lefteqn{
\left[
\begin{array}{cc} \Lambda & \ \\ \ & I \end{array}
\right]
{\left \langle \!\left[
\!\!\begin{array}{c} q(t) \\ p(t) \end{array}
\!\!\right] \!\right \rangle }_{exact,r} =
\left[ \begin{array}{cc} \cos {\Lambda t} & \sin {\Lambda t } \\
- \sin {\Lambda t } & \cos { \Lambda t }
\end{array} \right]
\left[ \begin{array}{cc} \Lambda^{-1} & \ \\ \ & I \end{array} \right]
\times} \\
& &
\left[ \begin{array}{cc} G & \ \\ \ & G \end{array} \right]
\left[ \begin{array}{cc}
(G^T\!\Lambda^{-2} G)^{-1} & \ \\
\ & (G^T\!G)^{-1}
\end{array} \right]
\left[ \!\! \begin{array}{c} G^T\!q \\ G^T\! p \end{array} \!\! \right].
\end{eqnarray*}
The index $r$ reminds us of the dimension.
Since $G^T = U [ I \ \cdots \ I ] \Gamma Y $
and $ Y \Lambda = \Lambda Y $, it follows that
\begin{displaymath}
{\langle \Lambda q(0) \rangle}_{exact,r} =
Y^\ast
\left[ \!\! \begin{array}{c} \Delta_{-r} \\ \vdots \\ \Delta_r \end{array}
\!\! \right]
\left[ \Delta_{-r}^2 + \cdots + \Delta_r^2 \right]^{-1}
\left[ \Delta_{-r} \ \cdots \ \Delta_r \right] Y \Lambda q ,
\end{displaymath}
where $\Delta_j = \Gamma_j \Lambda_j^{-1}$.
We get the formula for $\langle p(0) \rangle $ by replacing $\Lambda q $
by $p$ and $\Delta_j $ by $\Gamma_j $.
Next, let $P_r$ be the probability measure from Section 5 on
$ \Omega_r = {\Bbb R}^{2(2m+1)} $.
Since the random variables $a_i, \ b_i, \ \alpha_i , \ \beta_i $ are
independent, the measures $P_r$ are consistent and there is a probability space
$( \Omega, \cal{F} , P )$ such that
$P_r = P | \Omega_r $; see Billingsley \cite{Bil}, section 36.
We can now formulate

\begin{theorem}
\label{theorem2}
Let $n \geq 1$, and assume that
$(2n+1) \sigma^2 \geq 6(\nu+1) \log(2n+1) $ with $\nu \geq 0$.
Set $ \epsilon_r = 4(2n+1)^{-(\nu+1)[1+r+r^2(2n+1)]}$.
The limit of the exact method exists for almost all choices of the
random Fourier coefficients, and
\begin{displaymath}
\left| {\left \langle \!\!\left(
\!\!\!\begin{array}{c} u(x,t) \\ \pi (x,t) \end{array}
\!\!\!\right) \!\!\right \rangle }_{exact,r} -  \;\
{\left \langle \!\!\left(
\!\!\!\begin{array}{c} u(x,t) \\ \pi (x,t) \end{array}
\!\!\!\right) \!\!\right \rangle }_{exact,\infty} \right|_{A} <
\epsilon_r
\end{displaymath}
with probability greater than $1-\epsilon_r$.
\end{theorem}
{\bf Proof:}
Our proof is based on Borel-Cantelli, see \cite{Bil}, page 53.
Here is an outline.
Let $\psi_r (x,t,\omega) = {\langle {u(x,t) \choose \pi (x,t)}
\rangle}_{exact,r} $,
and define $\cal{A}_r= \{ \omega :| \psi_r-\psi_s|_{A} \geq \epsilon_r \}$.
Since $P(\cal{A}_r) < \epsilon_r $ and $\sum \epsilon_r < \infty$, it
follows that
$P( \cup_{s=0}^\infty \cap_{r=s}^\infty \cal{A}_r^c ) = 1 $.
The sequence $\{ \psi_r \}_{r=0}^{\infty}$
is therefore Cauchy for almost all $\omega \in \Omega $ and
converges to an element in $H^1\oplus H^0 $.
Instead of working with the random functions, we work with the
Fourier coefficients and imbed the smaller space into the larger space.
Let $r < s $, and set
\begin{eqnarray*}
B_1^T &=& [ \Delta_{-s} \cdots \Delta_s ] \\
B_2^T &=& [ \Delta_{-s} \cdots \Delta_{-r-1} \ 0 \cdots 0 \
\Delta_{r+1} \cdots \Delta_s ] \\
B_3^T &=& [ 0 \cdots 0 \ \Delta_{-r} \cdots \Delta_r \ 0 \cdots 0 ]
\end{eqnarray*}
Note that $B_i^T \! B_i $ are diagonal matrices of order $2n+1$.
We can now write
\begin{displaymath}
{\langle \Lambda q(0) \rangle}_{exact,s} -
{\langle \Lambda q(0) \rangle}_{exact,r} = b_1 + b_2 + b_3 ,
\end{displaymath}
where
\begin{eqnarray*}
b_1 &=& Y^\ast B_2 ( B_1^T \!B_1)^{-1} B_1^T Y \Lambda  q \\
b_2 &=& Y^\ast B_3 ( B_1^T \!B_1)^{-1} B_2^T Y \Lambda  q \\
b_3 &=&-Y^\ast B_3 ( B_1^T \!B_1)^{-1} B_2^T \!B_2 ( B_3^T \!B_3)^{-1} B_3^T
Y \Lambda q
\end{eqnarray*}
Using Chebyshev's inequality and Cauchy-Schwarz we see that
\begin{displaymath}
P(|\sum_{i=1}^3 b_i | > \epsilon ) \leq \epsilon^{-2}
E | \sum_{i=1}^3 b_i |^2 \leq 3 \epsilon^{-2} \sum_{i=1}^3 E|b_i|^2 .
\end{displaymath}
Since $\Lambda q $ are independent Gaussian random variables
with mean $0$ and variance $1$, and $ Y Y^\ast = I $ we get
\begin{displaymath}
E(b_1^T\! b_1) = \mathrm{tr} \ [Y^\ast B_1 (B_1^T \!B_1)^{-1} B_2^T \!
Y Y^\ast B_2 (B_1^T \!B_1)^{-1} B_1^T \!Y ] =
\mathrm{tr} \ [(B_1^T \!B_1)^{-1} B_2^T\! B_2 ].
\end{displaymath}
Now $\omega_k^2 / \omega_{k+j(2n+1)}^2 $ is less than $1$ if
$jk < 0 $ and less than $0.2$ if $jk \geq 0 $ and $j \neq 0 $.
Combining $ (B_1^T\!B_1)^{-1} < \Delta_0^{-2} $ with equation
(\ref{21})
and using $(2n+1)\sigma^2 \geq 6$ and $2n+1 \geq 3 $, we obtain
\begin{eqnarray*}
\lefteqn{
E|b_1|^2 \leq \sum_{k=-n}^n \sum_{j=r+1}^s
e^{-[(j^2-j)(2n+1)^2-j(2n+1)-2j(n-|k|)(2n+1) ] \sigma^2\!/2 }\ (1.2)} \\
& &
\leq e^{-r^2(2n+1)^2\sigma^2\!/2}
e^{-(r+1)(2n+1)\sigma^2\!/2}
\sum_{\ell =1}^\infty
(1.2) \, e^{-(\ell^2-\ell)3 \cdot 3 }
\sum_{k=-n}^n
e^{-(n-|k|)6} .
\end{eqnarray*}

Since $(2n+1)\sigma^2 \geq 6(\nu+1) \log(2n+1) $ and the product of the
two sums is less than $2.5$, we conclude that
\begin{displaymath}
E |b_1|^2 \leq 2.5 \, (2n+1)^{-3(\nu+1)[1+r+r^2(2n+1)]} = \epsilon' .
\end{displaymath}

By almost the same arguments, we get $E|b_2|^2 \leq \epsilon' $ for the
second term and for the third term we find that
$E|b_3|^2 \leq \mathrm{tr}\ [ ( \Delta_0^{-2} B_2^T\!B_2)^2 ]$.
Since $\Delta_0^{-2} B_2^T\!B_2$ is diagonal with all terms less than
$1$, it follows that $E|b_3|^2 \leq \epsilon' $.
The arguments for the $p$ terms are similar and by combining all estimates,
we obtain
\begin{equation}
\label{23}
{P} \left( | \psi_r - \psi_s |_A \geq (0.9) \epsilon_r \right)
\leq
(0.9 \epsilon_r)^{-2} \ 3 \cdot 2 \cdot 3 \cdot \epsilon'
< (0.9) \epsilon_r.
\end{equation}

Thus $P(\cal{A}_r) < \epsilon_r$.
Since $\sum \epsilon_r < \infty $, we conclude from Borel-Cantelli that
$P( \cap_{s=0}^{\infty} \cup_{r=s}^\infty \cal{A}_r ) = 0 $.
The sequence $\{ \psi_r \} $ is therefore a Cauchy sequence with probability
$1$ and $ \psi_r \rightarrow \psi_{\infty} $ in $H^1 \oplus H^0 $.
To estimate $ \psi_r - \psi_{\infty} $, we set
$\cal{B}_s = \cup_{r=s}^{\infty} \cal{A}_r $.
Since $\cal{B}_1 \supset \cal{B}_2 \supset \cdots $, there exists
an $s > r $ such that ${P} ( \cal{B}_s ) < (0.1) \epsilon_r $.
Let $\cal{A}_{rs} = \{ \omega : | \psi_r - \psi_s |_A < (0.9) \epsilon_r
\} $.
It follows from (\ref{23}) that
\begin{displaymath}
1-(0.9)\epsilon_r \leq
{P} ( \cal{A}_{rs} \cap \cal{B}_s^c ) +
{P} ( \cal{A}_{rs} \cap \cal{B}_s).
\end{displaymath}
The last term is less than $(0.1) \epsilon_r$, and for almost all
$\omega \in  \cal{A}_{rs} \cap \cal{B}_s^c $, we have
\begin{displaymath}
|\psi_r - \psi_{\infty} |_A \leq |\psi_r - \psi_s |_A +
|\psi_s - \psi_{s+1} |_A + \cdots <
(0.9)\epsilon_r + \sum_{j=s}^\infty \epsilon_j < \epsilon_r .
\end{displaymath}
Since $|(\psi_r - \psi_\infty )(x,t,\omega) |_A $ does not depend on
time, this completes the proof.

Suppose the constraints in (\ref{16}) are generated by
a smooth solution
${{u_0} \choose {\pi_0}}$ of (\ref{15}).
If $n\geq 1$ and $(2n+1)\sigma^2 \geq 6(\nu+1)
\log (2n+1)$ with $\nu \geq 0 $, we can show that
\begin{displaymath}
\left| {\left \langle \!\!\left(
\!\!\!\begin{array}{c} u(x,t) \\ \pi (x,t) \end{array}
\!\!\!\right) \!\!\right \rangle }_{exact,r} -  \;\
{\!\!\left(
\!\!\!\begin{array}{c} u_0(x,t) \\ \pi_0 (x,t) \end{array}
\!\!\!\right) } \right|_{A}
\end{displaymath}
\begin{displaymath}
\leq \frac{ 3 \sqrt{2.5} }{ (2n+1)^{(3/2)(\nu+1)}} \
\left| \left(\!\!\! \begin{array}{c} u_0 \\ \pi_0 \end{array}
\!\!\! \right) \right|_A
+
\frac{1}{(n+1)^s} \left| \ \partial_x^s ( I - Proj_n ) \left(\!\!\!
\begin{array}{c} u_0 \\ \pi_0 \end{array}
\!\!\! \right) \right|_A
\end{displaymath}
The method of optimal prediction can therefore also be used, in principle,
to solve the Klein-Gordon equation with smooth initial data.

\section*{Acknowledgments}
The author thanks Bradford Chin, Alexandre Chorin, Craig Evans,
William Kahan and Nicolai Reshetikhin for helpful discussions.
The work was supported in part by NSF under grant DMS-95-03483.

\bibliographystyle{plain}
\bibliography{ref}

\end{document}